\documentclass[a4paper,10pt]{article}
\usepackage{amsfonts}
\usepackage{amsmath}
\usepackage{amssymb}
\usepackage[latin1]{inputenc}
  \usepackage[T1]{fontenc}
  \usepackage{amsmath,amssymb}
  \usepackage[all]{xy}

\makeatletter

\@addtoreset{equation}{section} \makeatother

\newtheorem{thm}{Theorem}[section]
\newtheorem{cor}[thm]{Corollary}
\newtheorem{lem}[thm]{Lemma}
\newtheorem{prop}[thm]{Proposition}
\newtheorem{defn}{Definition}[section]
\newtheorem{rem}{Remark}[section]

\def\indlim{\mathop{\mathrm{ind\,lim}}}
\def\projlim{\mathop{\mathrm{proj\,lim}}}
\newenvironment{proof}[1][Proof]{\textbf{#1.} }
{\ \rule{0.5em}{0.5em}}

\begin{document}

\title{Quantum Gross Laplacian and Applications}
\author{{\bf Samah Horrigue}\; and \; {\bf Habib Ouerdiane}\\
Department of Mathematics\\ Faculty of sciences of Tunis\\
University of Tunis El-Manar\\ 1060 Tunis, Tunisia\\ e-mail:
samah.horrigue@fst.rnu.tn\\ e-mail: habib.ouerdiane@fst.rnu.tn }

\maketitle


\begin{abstract}
In this paper, we introduce and study a noncommutative extension of
the Gross Laplacian, called quantum Gross Laplacian. Then, applying
the quantum Gross Laplacian to the particular case where the
operator is the multiplication operator, we find a relation between
classical and quantum Gross Laplacian. As
application, we give explicit solution of linear quantum white noise
differential equation. In particular, we give a explicit solution of the quantum Gross heat equation.
\end{abstract}

\maketitle
\section{Introduction}\label{I}
The Gross Laplacian $\Delta_G$ was introduced by L. Gross in
\cite{g} in order to study differential equations in infinite
dimensional spaces. It has been shown that the solution of the
Cauchy problem
\begin{equation}\label{cauchyequation}
    \frac{\partial}{\partial t}U(t)=
    \frac{1}{2}\Delta_GU(t),\quad U(0)=\varphi
\end{equation}
is represented as an integral with respect to Gaussian measure, see
\cite{g} and \cite{pie70}. There exists many literature dedicated to
the Gross Laplacian with different points of view. We would like to
mention the white noise analysis approach, see \cite{CJ99, HKPS93,
HOS92, k'} and references therein. In \cite{coo} and \cite{bko},
using the fact that the Gross Laplacian is a convolution operator,
the authors applied Laplace transform techniques to solve the Cauchy
problem (\ref{cauchyequation}). Moreover, for a smooth initial condition the solution
is represented as an integral with respect to a Gaussian measure.

The main purposes of this work are the following: one is to define and study the generalized Gross Laplacian acting on operators denoted $\Delta_G^Q$. Another one is to solve linear quantum stochastic
differential equations. In particular, we give explicit solutions of the quantum Gross heat equations.

The paper is organized as follows. In section \ref{II}, we review
from \cite{joo} basic concepts, definitions and results essential to
know the space of test and generalized functions denoted
respectively $\mathcal{F}_ {(\theta_1, \theta_2)}(N'_1\times N_2')$
and $\mathcal{F}^\ast_ {(\theta_1, \theta_2)}(N'_1\times N_2')$. In
section \ref{III}, we introduce and define the
generalized Gross Laplacian $\Delta_G$ on the space of entire
functions with exponential growth of finite type in
 two infinite dimensional variables. Then, we prove in Theorem \ref{theorem2}, that
 the Gross Laplacian is a convolution operator on the test functions
 space $\mathcal{F}_ {(\theta_1, \theta_2)}(N'_1\times N_2')$, i.
 e., there exists a distribution $\mathbf{\mathcal{T}}\in\mathcal{F}^\ast_ {(\theta_1, \theta_2)}
 (N'_1\times N_2')$, such that
\begin{equation}\label{gross test}
\Delta_G\varphi=\mathbf{\mathcal{T}}\ast\varphi,\,\varphi\in
\mathcal{F}_{(\theta_1,\theta_2)}(N'_1\times N'_2).
\end{equation}
The relation (\ref{gross test}) permit to us to extend in natural way the Gross Laplacian applied to the
distributions spaces $\mathcal{F}^\ast_ {(\theta_1, \theta_2)}
 (N'_1\times N_2')$ as follows
\begin{equation}\label{gross gene}
\Delta_G\Phi=\mathbf{\mathcal{T}}\ast\Phi,\,\Phi\in
\mathcal{F}^\ast_{(\theta_1,\theta_2)}(N'_1\times N'_2).
\end{equation}
Using the Schwartz-Grothendieck Kernel Theorem and the definition
(\ref{gross gene}), we introduce and study, in section \ref{IV}, the noncommutative
extension of Gross Laplacian, denoted by $\Delta_G^Q$, acting on
$\mathcal{L}(\mathcal{F}_{(\theta_1,\theta_2)}(N'_1\times
N'_2),\mathcal{F}^\ast_{(\theta_1,\theta_2)}(N'
_1\times N'_2))$ the
space of continuous linear operators from
$\mathcal{F}_{(\theta_1,\theta_2)}(N'_1\times N'_2)$ into
$\mathcal{F}^\ast_{(\theta_1,\theta_2)}(N'_1\times N'_2)$. Then, we
establish an analytic characterization of the quantum Gross Laplacian (see
Proposition \ref{prop2}), i. e., for all
$\Xi\in\mathcal{L}(\mathcal{F}_{(\theta_1,\theta_2)}(N'_1\times
N'_2),\mathcal{F}^\ast_{(\theta_1,\theta_2)}(N'_1\times N'_2))$, we have
\begin{equation}\label{caracgross}
\sigma(\Xi)(\xi_1,\xi_2)=\left(\langle\xi_1,\xi_1\rangle_1+\langle
\xi_2,\xi_2\rangle_2\right)\sigma(\Xi)(\xi_1,\xi_2),\,(\xi_1,\xi_2)\in
N_1\times N_2,
\end{equation} where $\sigma(\Xi)$ denoted the symbol of the operator $\Xi$. In section \ref{V}, we study the action of the quantum Gross Laplacian to the
multiplication operator $\mathcal{M}_\Phi$ defined in
(\ref{defmulti}) where
$\Phi\in\mathcal{F}^\ast_{(\theta_1,\theta_2)}(N'_1\times N'_2)$, and we prove that:
\begin{equation}\label{0}
 \left(\Delta^Q_G\mathcal{M}_{\Phi}\right)e_0=
    \Delta_G\Phi,
\end{equation}
where $e_0$ is the vacuum vector. Therefore, the equality (\ref{0})
establish a relation between the classical and quantum Gross Laplacian.
In section \ref{VI},
 we give in Theorem \ref{t5}, the solution of the
following quantum stochastic differential equation
\begin{equation}\label{equdiff}
  (E)\left\{%
\begin{array}{ll}
\frac{d\Xi(t)}{dt}=Z(t)\ast\Xi(t)+\Theta(t),  \\
\Xi(0)=\Xi_0,\\
\end{array}%
\right.
\end{equation}
where $t\mapsto Z(t)$ and $t\mapsto\Theta(t)$ are  continuous operator valued process
defined on an interval containing the origin $I\subset\mathbb{R} $, i. e., \begin{eqnarray*}
                  Z:t\in I &\rightarrow&Z(t)\in\mathcal{L}(\mathcal{F}_ {\theta_1}(N_1')
    ,\mathcal{F}^\ast_ {
\theta_2}(N_2')), \\
           \Theta:t\in I &\rightarrow&\Theta(t)\in\mathcal{L}(\mathcal{F}_ {\theta_1}(N_1')
    ,\mathcal{F}^\ast_ {
\theta_2}(N_2'))
                \end{eqnarray*}
and the initial condition $\Xi_0\in\mathcal{L}(\mathcal{F}_ {\theta_1}(N_1')
    ,\mathcal{F}^\ast_ {
\theta_2}(N_2')).$
As an application of Theorem \ref{t5}, we give explicit solutions of the heat equation associated with the quantum Gross Laplacian.
\section{Preliminaries}\label{II}
For $i=1,2$, let $N_i$ be a complex nuclear Fr\'{e}chet space whose
topology is defined by a family of increasing Hilbertian norms
$\{|.|_{i, p},p\in\mathbb{N}\}$. For $p\in\mathbb{N}$, we denote by
$\left(N_i\right)_p$ the completion of $N_i$ with respect to the
norm $|.|_{i,p}$ and by $\left(N_i\right)_{i,-p}$ respectively
$N_i'$ the strong dual space of $(N)_p$ and $N$. Then, we obtain
\begin{equation}\label{1}
    N_i=\displaystyle \projlim_{p\rightarrow\infty}
    \left(N_i\right)_p\mbox{ and }N_i'=\displaystyle\indlim_
    {p\rightarrow\infty}\left(N_i\right)_{-p}.
\end{equation}
The spaces $N_i$ and $N_i'$ are respectively equipped with the
projective and inductive limit topology. For all $p\in\mathbb{N}$,
we denote by $|.|_{i,-p}$ the norm on $\left(N_i\right)_{-p}$ and by
$\langle.,.\rangle_i$ the $\mathbb{C}$-bilinear form on $N_i'\times
N_i.$

In the following, $H$ denote the direct Hilbertian sum of $(N_1)_0$ and
$(N_2)_0$, i. e., $H = (N_1)_0\oplus(N_2)_0$.

For $n\in\mathbb{N}$, we denote by $N_i^{\widehat{\otimes}n}$
the $n-$fold symmetric tensor product on $N_i$ equipped with the
$\pi-$topology and by $(N_i)^{\widehat{\otimes}n}_p$ the $n-$fold
symmetric Hilbertian tensor product on $(N_i)_p.$ We will preserve
the notation $|.|_{i,p}$ and
 $|.|_{i,-p}$ for the norms on $(N_i)^{\widehat{\otimes}n}_p$ and
 $(N_i)^{\widehat{\otimes}n}_{-p}$, respectively.

Let $\theta$ be a Young function, i. e.,
it is a continuous, convex and increasing function defined on $\mathbb{R}^+$ and satisfies the
 two conditions: $\theta(0)=0$ and
$\displaystyle\lim_{x\rightarrow+\infty}\frac{\theta(x)}{x}=+\infty$.
Obviously, the conjugate function $\theta^{\ast}$ of $\theta$
defined by
\begin{equation}\label{2}
    \forall
    x\geq0,\quad\theta^{\ast}(x):=\sup_{t\geq0}(tx-\theta(t)),
\end{equation}
is also a Young function.
For every $n\in\mathbb{N}$, let
\begin{equation}\label{5} \theta_n
=\displaystyle\inf_{r>0}\frac{e^{\theta(r)}}{r^n}.
\end{equation}
\subsection{Spaces of entire functions with growth condition}
Throughout the paper, we fix a pair of Young functions
$(\theta_1,\theta_2)$. For all pair of positive numbers $a_1,a_2>0$
and pair of integers $(p,q)\in\mathbb{N}\times\mathbb{N}$, we define
the space of all entire functions on
$\left(N_1\right)_{-p}\times\left( N_2\right )_{-q}$ with
$(\theta_1,\theta_2)-$exponential growth by
\[Exp(\left(N_1\right)_{-p}\times\left( N_2\right
)_{-q},(\theta_1,\theta_2),(a_1,a_2))= \{f\in\mathcal{H}(N_1\times
N_2);\|f\|_{(\theta_1,\theta_2),(a_1,a_2)}<\infty\},\] where
$\mathcal{H}(N_1\times N_2)$ is the space of all entire functions on
$N_1\times N_2$ and
\[\|f\|_{(\theta_1,\theta_2),(a_1,a_2)}= \sup\{|f(z_1,z_2)|
e^{-\theta_1(a_1|z_1|_{-p})-\theta_2(a_2|z_2|_{-q})},(z_1,z_2)\in
\left(N_1\right)_{-p}\times\left( N_2\right )_{-q}\}.\] So, the
space of all entire functions on $\left(N_1\right)_{-p}\times\left(
N_2\right )_{-q}$ with $(\theta_1,\theta_2)-$exponential growth of
minimal type is naturally defined by
 \begin{equation}\label{3}
    \mathcal{F}_{(\theta_1,
    \theta_2)}(N_1'\times N_2')=
        \displaystyle
    \projlim_{\stackrel{p,q\rightarrow\infty}{
    a_1,a_2\downarrow0}}Exp(\left(N_1\right)_{-p}\times\left( N_2\right
)_{-q},(\theta_1,\theta_2),(a_1,a_2)).
\end{equation}Similarly, the space of entire functions on $N_1\times N_2$
 with $(\theta_1,\theta_2)-$exponential growth of finite
type is defined by
\begin{equation}\label{4}
    \mathcal{G}_{(\theta_1,\theta_2)}(N_1\times N_2)=
    \displaystyle
    \indlim_{\stackrel{p,q\rightarrow\infty}{
    a_1,a_2\rightarrow0}}Exp(\left(N_1\right)_{p}\times\left( N_2\right
)_{q},(\theta_1,\theta_2),(a_1,a_2)).
\end{equation}

By definition, $\varphi\in
\mathcal{F}_{(\theta_1,\theta_2)}(N_1'\times N_2')$
 and $\Psi\in    \mathcal{G}_{(\theta_1,\theta_2)}(N_1\times N_2)$
  admit the Taylor expansions:

  \begin{equation}\label{taylor1}
     \varphi(x,y)=\displaystyle\sum_{n,m\in\mathbb{N}}\langle
      x^{\otimes n}\otimes y^{\otimes m},\varphi_{n,m}
      \rangle,\,(x,y)\in N_1'\times N_2',
  \end{equation}
and\begin{equation}\label{taylor2}
    \Psi(\xi,\eta)=\displaystyle\sum_{n,m\in\mathbb{N}}\langle
\Psi_{n,m},      \xi^{\otimes n}\otimes \eta^{\otimes m}
      \rangle,\,(\xi,\eta)\in N_1\times N_2,
\end{equation}
where for all $n,m\in\mathbb{N}$,
we have $\varphi_{n,m}
      \in N_1^{\hat{\otimes} n}\otimes N_2^{\hat{\otimes} m}$, $ \Psi_{n,m}
      \in \left(N_1^{\hat{\otimes} n}\right)'\otimes \left
      (N_2^{\hat{\otimes} m}\right)'$ and we used the common
       symbol $\langle\,.\,,\,.\,\rangle
$ for the canonical $\mathbb{C}-$bilinear form on
 $\left(N_1^{\otimes n}\times N_2^{\otimes m}\right)'\times
N_1^{\otimes n}\times N_2^{\otimes m}$.
So, we identify in the next all test function $\varphi\in\mathcal{F}
_{(\theta_1,\theta_2)}(N'_1\times N'_2)$ (resp. all generalized function
$\Psi\in\mathcal{G} _{(\theta_1,\theta_2)}(N_1\times N_2)$) by their coefficients of its
Taylor series expansion at the origin $ (\varphi_{n,m})
_{n,m\in\mathbb{N}}$
(resp. $(\Psi_{n,m}) _{n,m\in\mathbb{N}}$).

Denote by $
\mathcal{F}^\ast_{(\theta_1,\theta_2)}(N_1'\times N_2')$
 the topological dual of
$\mathcal{F}_{(\theta_1,\theta_2)}(N_1'\times N_2')$,
 called the space of distributions on $N_1'\times
N_2'$.

For a fixed $(\xi,\eta)\in N_1\times N_2$, the exponential function
$e_{(\xi,\eta)}\in\mathcal{F}_{(\theta_1,\theta_2)}(N_1'\times N_2')$
  is defined by
\[e_{(\xi,\eta)}(z,t)=\exp\{\langle z,\xi\rangle_1+\langle
t,\eta\rangle_2\} ,\,(z,t)\in N_1'\times N_2'.\]
Then for every $\Phi\in\mathcal{F}^\ast_{(\theta_1,\theta_2)}
(N_1'\times N_2')$, the Laplace transform $L$ of $\Phi$ is defined by
\begin{equation}\label{2'}
L\Phi(\xi,\eta)=\widehat{\Phi}(\xi,\eta)=\langle\langle\Phi
,e_{(\xi,\eta)}\rangle\rangle.
\end{equation}
\begin{thm}\label{t3}\cite{joo}
For $i=1,2$, let $N_i$ be complex nuclear Fr\'{e}chet space and
$\theta_i$ a Young function. Then, the Laplace transform $L$ is a
topological isomorphism:
\begin{equation}\label{laplace}
L:\mathcal{F}^\ast_{(\theta_1,\theta_2)}(N_1'\times N_2')
\rightarrow\mathcal{G}_{(\theta_1,\theta_2)}(N_1\times
N_2) .
\end{equation}

\end{thm}
\begin{rem}\label{rem1}
In the particular case where $\theta_1=\theta_2
 =\theta$, $N_1=N$ and $N_2=\{0\}.$ We obtain the following identification
have\[\mathcal{F}_ {(\theta,\theta)}(N'\times \{0\})=\mathcal{F}_
{\theta}(N')\] and therefore \[\mathcal{F}^\ast_ {(\theta,\theta)}(N'\times
\{0\}) =\mathcal{F}^\ast_ {\theta}(N').\]
So the space $\mathcal{F}_ {(\theta_1,\theta_2)}(N'\times N_2')$ can be considered as a generalization of the space $\mathcal{F}_ {\theta}(N')$ studied in \cite{ghor}.
\end{rem}
\subsection{Convolution operators}
Let $(z,t)\in
    N_1'\times N_2'$, the translation operator denoted $\tau_{(z,t)}$ is a linear continuous operator
    $\tau_{(z,t)}$ from $\mathcal{F}
_{(\theta_1,\theta_2)}(N_1'\times N_2')$ into itself is defined by
\begin{equation}\label{21}
\tau_{(z,t)}\varphi(x,y)=\varphi(x+z,y+t),\,\, (x,y)\in
    N_1'\times N_2'.
\end{equation}

The convolution product of the distribution
$\Phi\in\mathcal{F}^\ast_{(\theta_1,\theta_2)}(N_1'\times N_2')$
with a test function
$\varphi\in\mathcal{F}_{(\theta_1,\theta_2)}(N_1'\times N_2')$ is
defined as follows
\begin{equation}\label{20}
    \Phi\ast\varphi(z,t)=\langle\langle\Phi,\tau_{(z,t)}
    \varphi\rangle\rangle,\,\,(z,t)\in
    N_1'\times N_2'.
\end{equation}

Hence, the convolution product of two distribution
$\Psi_1,\Psi_2\in\mathcal{F}^\ast_{(\theta_1,\theta_2)}(N_1'\times
N_2')$ is given by
\begin{equation}\label{22}
    \widehat{\Psi_1\ast\Psi_2}=\widehat{\Psi_1}\widehat{\Psi_2}
    ,\,\Psi_1,\Psi_2\in\mathcal{F}^\ast_{(\theta_1,\theta_2)}(N_1'\times
N_2').
\end{equation}So, by formula (\ref{22}), the convolution product of distribution is commutative and associative.

We denote by $\mathcal{L}(\mathcal{F} _ {\theta_1}(N_1')
    ,\mathcal{F}^\ast_ {\theta_2}(N_2'))$ the space of continuous linear operators from $\mathcal{F}
_{\theta_1}(N_1')$ into $\mathcal{F}^\ast _{\theta_2}(N_2')$ endowed
with the bounded convergence topology.

A convolution operator on the test
space $\mathcal{F} _{(\theta_1,\theta_2)}(N_1'\times N_2')$ is defined as a
continuous linear operator from
$\mathcal{F}_{(\theta_1,\theta_2)}(N_1'\times N_2')$ into itself
which commutes with translation operators. Then, $T$ is a
convolution operator on $\mathcal{F}
_{(\theta_1,\theta_2)}(N_1'\times N_2')$ if and only if there exists
a distribution $\Phi_T\in\mathcal{F}
    ^\ast_{(\theta_1,\theta_2)}(N_1'\times
N_2')$ such that
\[T(\varphi)=T_\Phi(\varphi)=\Phi_T\ast\varphi,\,\varphi\in\mathcal{F}
_{(\theta_1,\theta_2)}(N_1'\times N_2').\]
Note that the convolution product of two distributions $\Psi_1,\Psi_2\in\mathcal{F}
    ^\ast_{(\theta_1,\theta_2)}(N_1'\times
N_2')$ is also given by
\[\langle\!\langle\Psi_1\ast\Psi_2,\varphi\rangle\!\rangle
=\langle\!\langle\Psi_1,\Psi_2\ast\varphi\rangle\!\rangle,\,
\varphi\in\mathcal{F}
_{(\theta_1,\theta_2)}(N_1'\times
N_2').\]

It follows that the notion of
convolution operator $T_\Phi$ can be extended to the distribution space
$\mathcal{F}
    ^\ast_{(\theta_1,\theta_2)}(N_1'\times
N_2')$ as follows: $\mathbf{T}$ is a convolution operator on $\mathcal{F}
    ^\ast_{(\theta_1,\theta_2)}(N_1'\times
N_2')$ if and only if there exist $\Phi_T\in\mathcal{F}
    ^\ast_{(\theta_1,\theta_2)}(N_1'\times
N_2')$ such that
\[ \mathbf{T}(\Psi)= \mathbf{T}_\Phi(\Psi)=\Phi_T\ast\Psi,\,\Psi\in\mathcal{F}^\ast
_{(\theta_1,\theta_2)}(N_1'\times N_2').\]
We remark that this extension of convolution operator $\mathbf{T}$ coincide with the adjoint of $T_\Phi$ denoted $T^\ast_\Phi$. In fact, for all $\Psi\in\mathcal{F}
    ^\ast_{(\theta_1,\theta_2)}(N_1'\times
N_2')$ and $\varphi\in\mathcal{F}
_{(\theta_1,\theta_2)}(N_1'\times
N_2')$, we have:\begin{eqnarray*}
\langle\!\langle T^\ast_\Phi\Psi
,\varphi\rangle\!\rangle&=:& \langle\!\langle \Psi
,T_\Phi\varphi\rangle\!\rangle \\
                   &=& \langle\!\langle \Psi
,\Phi\ast\varphi\rangle\!\rangle
 \\
&=& \langle\!\langle \mathbf{T} _\Phi
\Psi,\varphi\rangle\!\rangle.
                \end{eqnarray*}

Let $i=1,2$. Recall now the notion of the right contraction of order
$k\in\mathbb{N}$:

 For $n,m\in\mathbb{N}-\{0\}$ and
$0\leq k\leq m\wedge n$, we denote by $\langle\,.,\,.\rangle_{i,(k)}$
the bilinear map from $N_i'^{\hat{\otimes}m}\times
N_i^{\hat{\otimes}n}$ into
$N_i'^{\hat{\otimes}m-k}\widehat{\otimes}N_i^{\hat{\otimes}n-k}$
defined by
\[\langle x^{\otimes m},y^{
\otimes n}\rangle_{i,(k)}:=\langle x,y\rangle_i^k x^{\otimes(m-k)}\otimes
y^{ \otimes(n-k)},\,x\in N_i', y\in N_i.\] The bilinear map
$\langle\,.,\,.\rangle_{i,(k)}$ is continuous. Then using the density
of the vector space generated by $\{x^{\otimes m},\, x\in N_i'\}$ in
$N_i'^{\hat{\otimes}m}$ and the vector space generated by $\{y^{
\otimes n},\,y\in N_i\}$ in $N_i^{\hat{\otimes}n}$, we can extend
$\langle\,.,\,.\rangle_{i,(k)}$ to $N_i'^{\hat{\otimes}m}\times
N_i^{\hat{\otimes}n}$. For all $\Phi_n\in N_i'^{\hat{\otimes}m}$ and
$\varphi_n\in N_i^{\hat{\otimes}n}$,
$\langle\Phi_n,\varphi_n\rangle_{i,(k)}$ is called the right contraction
of $\Phi_n$ and $\varphi_n$ of order $k$.

 It follows that for
$(\varphi_1,\varphi_2)\in N_1^{\widehat{\otimes}(n+i)}\times
N_2^{\widehat{\otimes}(m+j)}$ and $(\psi_1,\psi_2)\in
N_1^{\widehat{\otimes}(n+k)}\times N_2^{\widehat{\otimes}(m+l)}$, we
define the following generalized contraction as follows:
\[\langle \varphi_1\otimes\varphi_2,\psi_1\otimes \psi_2\rangle_{n,m}=
\langle \varphi_1,\psi_1\rangle_{1,(n)}\langle \varphi_2,\psi_2\rangle_{2,(m)}.\] So, for
$i=1,2$, $\Phi_i\in (N'_{i})^{\widehat{\otimes}(p_i+m_i)}$, $\varphi_i\in
N_i^{\widehat{\otimes} (p_i+n_i)}$ and $\psi_i\in
N_i^{\widehat{\otimes}(m_i+n_i)}$, it holds:
\[\langle \langle \Phi_1\otimes\Phi_2,\varphi_1\otimes \varphi_2
\rangle_{p_1,p_2},\psi_1\otimes\psi_2\rangle=\langle \Phi_1\otimes
\Phi_2,\langle \varphi_1\otimes\varphi_2, \psi_1\otimes
\psi_2\rangle_{n_1,n_2}\rangle.\]By an easy calculation like in
\cite{co}\,, we obtain the following lemma.
\begin{lem}\label{lem2}
For all test function
$\varphi=(\varphi_{n,m})_{n,m\in\mathbb{N}}\in\mathcal{F}_{(\theta_1,\theta_2)}(N_1'\times
N_2')$ and generalized function
$\Phi=(\Phi_{n,m})_{n,m\in\mathbb{N}}\in\mathcal{F}^\ast_{(\theta_1,\theta_2)}(N_1'\times
N_2')$, we have
\begin{align}\label{taylorconvolution}
    & \Phi\ast\varphi(x,y)\nonumber\\&
     =\displaystyle\sum_{k,l\in\mathbb{N}}
    \frac{(n+k)!}{k!}    \frac{(m+l)!}{l!}    \left\langle
    x^{\otimes k}\otimes y^{\otimes l}
    ,\langle\displaystyle\sum_{n,m\in\mathbb{N}}
\Phi_{n,m}
    ,
    \varphi_{n+k,m+l}\rangle_{n,m}\right\rangle,\,\end{align}
for all $(x
    ,y)\in N'_1\times N'_2$.
\end{lem}

\section{Generalized Gross Laplacian}\label{III}
Let $F\in C^2(N)$. Then for each $\xi\in N$ there exist $F^\prime(\xi)\in N'$ and $F^{\prime\prime}(\xi)\in(N\otimes N)'$ such that
\begin{equation}\label{devtay1}
    F(\xi+\eta)=F(\xi)+\langle F^{\prime}(\xi),\eta\rangle+\frac{1}{2}\langle F^{\prime\prime}(\xi),\eta\otimes\eta\rangle+o(|\eta|^2_p)
,\,\eta\in N,
\end{equation}
for some $p\in\mathbb{N}$. Moreover, both maps $\xi\mapsto F^{\prime}(\xi)\in N'$ and $\xi\mapsto F^{\prime\prime}(\xi)\in (N\otimes N)'$ are continuous. For notation simplicity, taking into account the canonical isomorphism $(N\otimes N)'\simeq\mathcal{L}(N,N')$, which follows from the kernel theorem for a nuclear space, we write $\langle F^{\prime\prime}(\xi),\eta\otimes\eta\rangle=\langle F^{\prime\prime}(\xi)\eta,\eta\rangle= F^{\prime\prime}(\xi)(\eta,\eta).$

Let  $\tau$ be the trace operator defined on $N^{\otimes
2}$ by the formula
\begin{equation}\label{14}
\langle\tau,\xi\otimes\eta\rangle=\langle\xi,\eta\rangle,
\forall\xi,\eta\in N.
\end{equation}

Using the definition of the Gross Laplacian
 $\Delta_G$ given in \cite{g}, the authors studied in
 \cite{co}, the action of $\Delta_G$
on $\mathcal{F} _{ \theta}(N')$.
Then, they prove that $\Delta_G$ is in fact a convolution
 operator i. e., for all $\varphi\in\mathcal{F}
 _{ \theta}(N')$, we have
\begin{eqnarray}
\Delta_G\varphi(x) &:=& trace_HD^2\varphi(x)\nonumber\\
  &=&\displaystyle\sum_{n\in\mathbb{N}}(n+2)(n+1)
    \langle x^{\otimes^n},
    \langle\tau,\varphi_{n+2}
    \rangle\rangle \label{19"}\\
&=&  \mathcal{T}\ast\varphi(x),\,\,x\in N'\label{18},
\end{eqnarray}
where the Taylor expansion of the distribution $\mathcal{T}$ is
given by
\begin{equation}\label{24}
\mathcal{T}_{n}=\left\{%
\begin{array}{ll}
        0,\,n\neq2\\
\tau,\,n=2. \\
\end{array}%
\right.
\end{equation}

Then in the paper \cite{bko}\,, the authors extended the action of the
Gross Laplacian on the test function space given in (\ref{18}) to the distribution
space as follows
\begin{equation}\label{19'"}
\Delta_G\Phi=\mathcal{T}\ast\Phi,\,\Phi\in \mathcal{F}^\ast _{
\theta}(N').
\end{equation}
In this section, we define first the action of the Gross Laplacian on the test function space of
two infinite dimensional variables $\mathcal{F} _{(
\theta_1,\theta_2)}(N_1'\times N_2')$.

For $i=1,2,$ let $\{e^{i}_{j}\}_{j\in\mathbb{N}}$ be a complete
orthonormal basis of $(N_i)_0$ such that $e^{i}_{j}\in N_i $. In all
the remainder of this part, we denote by $\tau_i$ the trace operator on
$N_i^{\otimes2}$ defined in (\ref{14}).
\begin{lem} \cite{ob}\label{lem1}  For all $i=1,2$ and $w_i \in N^{\otimes2}_i$, it holds that
\[\langle\tau_i,w_i\rangle_i
=\displaystyle\sum_{n\in\mathbb{N}} \langle (e_n^i)^{\otimes^2},w_i
    \rangle_i.\]
\end{lem}
We recall that all function in two variables $(\xi_1,\xi_2)\in N_1\times N_2$ is identified in an obvious manner with a single-variable function on the direct sum $N=N_1\oplus N_2$, which is again a countably Hilbert nuclear space. Then, for $F\in C^2(N)$, (\ref{devtay1}) is written, for some $p,q\in\mathbb{N}$, in the following form:

\begin{eqnarray}
  F(\xi_1+\eta_1,\xi_2+\eta_2)&=&
 F(\xi_1,\xi_2)
+\displaystyle\sum_{i=1}^2\langle F^{\prime}_i(\xi_1,\xi_2),\eta_i\rangle+
\frac{1}{2}\displaystyle\sum_{i,j=1}^2\langle F^{\prime\prime}_{i,j}(\xi_1,\xi_2)\eta_i,\eta_j\rangle \nonumber \\
        &+& o(|\eta_1|_p^2
+|\eta_2|_q^2),\label{devtay2}
     \end{eqnarray}
where $F^{\prime}_i(\xi_1,\xi_2)\in N'_i$,
 $F^{\prime\prime}_{i,j}(\xi_1,\xi_2)\in \mathcal{L}(N_i,N_j')$ for $i,j=1,2$ and the error term satisfies $$\displaystyle\lim_{t\rightarrow0}\frac{o(t(|\eta_1|_p^2
+|\eta_2|_q^2))}{t^2}=0.$$

For $i=1,2$, identifying $(N_i)_0$ with a subspace of
$H=(N_1)_0\oplus(N_2)_0$ in the canonical manner, we write $e_n^1$ and $e_n^2$ for $e_n^1\oplus0$ and $0\oplus e_n^2$, respectively, for all $n\in\mathbb{N}$.

\begin{thm}\label{theorem2}
For any test function $\varphi(x,y)=\displaystyle\sum_{n,m\in\mathbb
{N}}\langle x^{\otimes n}\otimes y^{\otimes
m},\varphi_{n,m}\rangle$ in $\mathcal{F}_{(\theta_1,\theta_2)}(N_1'\times
N_2')$, the Gross Laplacian is given by:
\begin{eqnarray}
  \Delta_G\varphi (x,y)   &:=& trace_HD^2\varphi(x,y)  \nonumber \\
   &=&\displaystyle\sum_{n,m\in\mathbb{N}}
    \langle x^{\otimes^n}\otimes y^{\otimes^m},
    \langle(n+2)(n+1)\tau_1,\varphi_{n+2,m}
    \rangle\!\rangle \nonumber  \\&+&\displaystyle\sum_{n,m\in\mathbb{N}}
    \langle x^{\otimes^n}\otimes y^{\otimes^m},
    \langle(m+2)(m+1)\tau_2,\varphi_{n,m+2}
    \rangle\!\rangle\label{19}.
\end{eqnarray}
Moreover, the Gross Laplacian $\Delta_G$ is a convolution operator on $\mathcal{F}_{(\theta_1,\theta_2)}(N_1'\times
N_2')$ into itself, i. e., there exists a distribution  $\mathbf{\mathcal{T}}=\left(\mathbf{\mathcal{T}}_{n,m}\right)_{n,m\in
\mathbb{N}}$ such that
\begin{equation} \Delta_G(\varphi)=\mathbf{\mathcal{T}}\ast\varphi,\,\varphi
\in\mathcal{F}_{(\theta_1,\theta_2)}(N_1'\times
N_2'),\label{21'}
\end{equation}
where the Taylor expansion of the distribution $\mathbf{\mathcal{T}}$ is
given by
\begin{equation}\label{24}
\mathbf{\mathcal{T}}_{n,m}=\left\{%
\begin{array}{lll}
       \tau_1,\,n=2,\,m=0,\\
\tau_2,\,n=0,\,m=2, \\
0,\,\mbox{otherwise.}
\end{array}%
\right.
\end{equation}
\end{thm}
\textbf{Proof.}
Let $\varphi(x_1,x_2)=\displaystyle\sum_{n,m\in
  \mathbb{N}} \langle x_1^{\otimes n}
  \otimes x_2^{\otimes m},\varphi_{n,m}\rangle$ be a test function in
$\mathcal{F}_ {(\theta_1,\theta_2)}(N_1'\times N_2')$. So, for all
$(x_1,x_2),(y_1,y_2)\in N_1'\times N_2'$, the Taylor expansion of the function $\varphi(x_1+y_1,x_2+y_2)$ is given by

\begin{align}\label{7}
    &\varphi(x_1+y_1,x_2+y_2)
    =\displaystyle\sum_{n,m\in
  \mathbb{N}} \langle (x_1+y_1)^{\otimes^n}
  \otimes (x_2+y_2)
  ^{\otimes^m},\varphi_{n,m}\rangle\nonumber
  \\&\qquad\displaystyle=\sum_{n,m,i,j\in
  \mathbb{N}} \left(%
\begin{array}{c}
  n+i \\
  i \\
\end{array}%
\right)\left(%
\begin{array}{c}
  m+j \\
  j \\
\end{array}%
\right)\langle x_1^
{\otimes^n}\otimes
x_2
  ^{\otimes^m},\langle y_1^{\otimes
  i}\otimes y_2^{\otimes
  j},\varphi_{n+i,m+j}\rangle_{i,j}\rangle\nonumber\\&\,=
  \varphi(x_1,x_2)+
  \displaystyle\sum_{n,m\in
  \mathbb{N}}\langle x_1^{\otimes^n}\otimes
   x_2
  ^{\otimes^m},\langle y_1,(n+1)\varphi_{n+1,m}\rangle_{1,0}+
\langle y_2,(m+1)\varphi_{n,m+1}\rangle_{0,1}\rangle\nonumber\\&+
\displaystyle\sum_{i,j=1}^2
  \displaystyle\sum_{n,m\in\mathbb{N}} \left(%
\begin{array}{c}
  n+i \\
  i \\
\end{array}%
\right)\left(%
\begin{array}{c}
  m+j \\
  j \\
\end{array}%
\right)
\langle
x_1^{\otimes^n}\otimes x_2
  ^{\otimes^m},\langle y_1^{\otimes
  i}\otimes y_2^{\otimes
  j},\varphi_{n+i,m+j}\rangle_{i,j}\rangle+\varepsilon
  (y_1,y_2),
  \end{align}where for each $(x_1,x_2)\in N_1'\times N_2'$, $\varepsilon$ is given by
  \[\varepsilon
  (y_1,y_2)=\displaystyle\sum_{\stackrel{n,mi,j\in
  \mathbb{N}}{i,j\geq3}} \left(%
\begin{array}{c}
  n+i \\
  i \\
\end{array}%
\right)\left(%
\begin{array}{c}
  m+j \\
  j \\
\end{array}%
\right)
\langle
x_1^{\otimes^n}\otimes x_2
  ^{\otimes^m},\langle y_1^{\otimes
  i}\otimes y_2^{\otimes
  j},\varphi_{n+i,m+j}\rangle_{i,j}\rangle.
\]It is easy to see that the function $\varepsilon$ satisfy \[\displaystyle\lim_{t
\rightarrow0}\frac{\varepsilon\left(t(y_1,  y_2)\right)}{t}=0.\]
Identifying the equality (\ref{7}) to the development (\ref{devtay2}), the second
derivative of $\varphi$ at $(x_1 ,x_2)$ in the direction $(y_1,y_2)\in H$ is given by
\[
  D^2\varphi(x_1,x_2)=\displaystyle\sum_{i,j=1}^2
  \displaystyle\sum_{n,m\in\mathbb{N}} \left(%
\begin{array}{c}
  n+i \\
  i \\
\end{array}%
\right)\left(%
\begin{array}{c}
  m+j \\
  j \\
\end{array}%
\right)
\langle
x_1^{\otimes^n}\otimes x_2
  ^{\otimes^m},\langle y_1^{\otimes
  i}\otimes y_2^{\otimes
  j},\varphi_{n+i,m+j}\rangle_{i,j}\rangle.\]
Using Lemma \ref{lem1}, we have for all $(x,y)\in N_1'\times
N'_2$:

\begin{align}
&trace_{H}D^2\varphi(x,y)=\displaystyle
\sum_{i,j\in\mathbb{N}}\langle D^2\varphi(x,y)(e^1_i,e^2_j),e^1_i\otimes e^2_j\rangle\nonumber\\
   &\,=\displaystyle
\sum_{n,m,i,j\in\mathbb{N}}
    \langle x^{\otimes^n}\otimes y^{\otimes^m},
 (n+2)(n+1) \langle (e^1_i)^{\otimes2},\varphi_{n+2,m}
    \rangle_{2,0}+ (m+2)(m+1)\langle (e^2_j)^{\otimes 2},\varphi_{n,m+2}
    \rangle_{0,2}\rangle\nonumber\\
&\,=\displaystyle\sum_{n,m\in\mathbb{N}}
    \langle x^{\otimes^n}\otimes y^{\otimes^m},
    \langle((n+2)(n+1)\tau_1,\varphi_{n+2,m}
    \rangle_{2,0}+ (m+2)(m+1)\tau _2,\varphi_{n,m+2}
    \rangle_{0,2})\rangle \label{6}.
\end{align}
In the other hand, using Lemma \ref{lem2} for the distribution $\mathbf{\mathcal{T}}\in\mathcal{F}^\ast_
{(\theta_1,\theta_2)}(N_1'\times N_2')$ where the Taylor expansion is given in (\ref{24}) and the equality (\ref{6}), we obtain
\begin{align*}
    & \mathbf{\mathcal{T}}\ast\varphi(x,y)  =\displaystyle\sum_{n,m\in\mathbb{N}}
        \langle x^{\otimes^n}\otimes y^{\otimes^m},\langle
        \mathbf{\mathcal{T}}_{2,0}        ,\varphi_{n+2,m}\rangle_{2,0}+        \mathbf{\mathcal{T}}_{0,2}
        ,\varphi_{n,m+2}\rangle_{0,2}\rangle\\&
    = \displaystyle\sum_{n,m\in\mathbb{N}}
    \langle x^{\otimes^n}\otimes y^{\otimes^m},
    \langle((n+2)(n+1)\tau_1,\varphi_{n+2,m}
    \rangle_{2,0}+ (m+2)(m+1) \tau _2),\varphi_{n,m+2}
    \rangle_{0,2}\rangle \\
&=\Delta_G\varphi(x,y),\,(x,y)\in N_1'\times N_2',\,\varphi\in\mathcal{F}_
{(\theta_1,\theta_2)}(N_1'\times N_2').
\end{align*}
Since the Gross Laplacian $\Delta_G$ is a convolution operator then $$\Delta_G\in\mathcal{L}(\mathcal{F}_{(\theta_1
      ,\theta_2)}(N'_1\times N'_2),\mathcal{F}_{(\theta_1
      ,\theta_2)}(N'_1\times N'_2)).$$
\begin{rem}
Let $\varphi\in\mathcal{F}_{(\theta_1
      ,\theta_2)}(N'_1\times N'_2)$. Then, the Taylor expansion (\ref{19}) of $\Delta_G$, can be written as follows
\begin{equation}\label{18'}
    \Delta_G\varphi=\Delta_G^1\varphi+
    \Delta_G^2\varphi,
\end{equation}
where for $i=1,2$, $\Delta_G^i\varphi$ is the Gross Laplacian with respect to the variable $\xi_i\in N'_i$ given by
\[\Delta_G^1\varphi(x_1,x_2)=\displaystyle
\sum_{n,m\in\mathbb{N}} \langle x_1^{\otimes^n}\otimes x_2^{\otimes^m},
    \langle(n+2)(n+1)\tau_1,\varphi_{n+2,m}
    \rangle\rangle\]respectively,
    \[\Delta_G^2\varphi(x_1,x_2)=\displaystyle
\sum_{n,m\in\mathbb{N}} \langle x_1^{\otimes^n}\otimes x_2^{\otimes^m},
    \langle(m+2)(m+1)\tau_1,\varphi_{n,m+2}
    \rangle\rangle.\]
\end{rem}Using the previous remark, it is easy to prove the following result.
\begin{prop}\label{prop2}
For $\varphi=f\otimes g\in\mathcal{F}_{\theta_1}(N'_1)\otimes
\mathcal{F}_{\theta_2}(N'_2)\subset\mathcal{F}_{(\theta_1,\theta_2)}
    (N'_1\times N'_2)$, we obtain
    \[\Delta_G\varphi=\Delta_G(f\otimes g)=\Delta_G(f)\otimes g+
    f\otimes\Delta_G(g).\]Therefore, on the subspace $\mathcal{F}_{\theta_1}(N'_1)\otimes
\mathcal{F}_{\theta_2}(N'_2)$ of $\mathcal{F}_{(\theta_1,\theta_2)}(N'_1\times
N'_2)$, the Gross Laplacian $\Delta_G$ is given by
\[\Delta_G=(\Delta_G)_{|\mathcal{F}_{\theta_1}(N'_1)}\otimes I
_{2}+I_{1}\otimes(\Delta_G)_{|\mathcal{F}_{\theta_2}(N'_2)} ,\] where for $i=1,2$, $(\Delta_G)_{|\mathcal{F}_{\theta_i}(N'_i)}$ (respectively $I_i$) is the classical Gross Laplacian (respectively the identity operator) acting on $\mathcal{F}_{\theta_i}(N'_i)$.
\end{prop}
Then, we obtain the following corollary:
\begin{cor}\label{cor2}
For all $(\xi_1,\xi_2)\in N_1\times N_2$, we have
\begin{equation}\label{grossexponentiel}
    \Delta_G(e_{(\xi_1,\xi_2)})=    \Delta_G(e_{\xi_1}\otimes e_{\xi_2})=(\langle\xi_1,\xi_1\rangle_1+
    \langle\xi_2,\xi_2\rangle_2)e_{(\xi_1,\xi_2)}.
\end{equation}
\end{cor}
\begin{proof}To prove the corollary, it is sufficient to see that for $i=1,2$, we have:
\begin{eqnarray*}
  (\Delta_G)_{|\mathcal{F}_{\theta_i}(N'_i)}
  (e_{\xi_i})(x_i) &=&\displaystyle\sum_{n
\in\mathbb{N}} \langle x_i^{\otimes n},\langle\tau_i,\xi_i^{\otimes n+2}\rangle_2\rangle \nonumber \\
   &=& \langle\xi_i,\xi_i\rangle_i\displaystyle
   \sum_{n
\in\mathbb{N}}\langle x_i^{\otimes n},,\xi_i^{\otimes n}\rangle \nonumber\\
   &=&\langle\xi_i,\xi_i\rangle_i e_{\xi_i}(x_i) , \nonumber
\end{eqnarray*}for all $\xi_i\in N_i$ and $x_i\in N'_i$.
\end{proof}

For $i=1,2$, we assume that the Young function
$\theta_i$ satisfies the following condition:
\begin{equation}\label{c1}
    \displaystyle\limsup_{x\rightarrow\infty}
    \frac{\theta_i(x)}{x^2}<+\infty,
\end{equation}
we obtain the following Gel'fand triple (see \cite{ghor})

\begin{equation}\label{triplet1}
\mathcal{F}_{\theta_i}(N_i')\rightarrow L^2(X_i',\gamma_i)
\rightarrow\mathcal{F}^\ast_{\theta_i}(N_i'),
\end{equation}
\begin{equation}\label{triplet2}
    \mathcal{F}_{(\theta_1,\theta_2)}(N_1'\times
    N_2')\rightarrow L^2(X_1'\times
    X_2',\gamma_1\otimes\gamma_2)
\rightarrow\mathcal{F}^\ast_{(\theta_1,\theta_2)}(N_1'\times
    N_2'),
\end{equation}

 where $\gamma_i$ is the Gaussian measure on the real Fr\'{e}chet nuclear space $X_i'$ whose complexification is $N_i$ defined via the
Bochner-Minlos theorem \cite{HKPS93} by its characteristic
function:
\begin{equation}\label{gaussian measure}
\int_{X'_i}e^{i\langle
x,\xi_i\rangle_i}d\gamma(x)=e^{-\frac{1}{2}|\xi_i|_0^2}, \xi_i\in X_i.
\end{equation}

Using the Gel'fand triplet (\ref{triplet2}) and the fact that the Gross Laplacian is a convolution product, we extend the Gross Laplacian on $\mathcal{F}^\ast_
{(\theta_1,\theta_2)}(N_1'\times N_2')$ as follows:
\begin{defn}The generalized Gross Laplacian acting on the distributions space $\mathcal{F}^\ast_
{(\theta_1,\theta_2)}(N_1'\times N_2')$ is defined by
\begin{equation}\label{23}
    \Delta_G(\Psi)=\mathbf{\mathcal{T}}\ast\Psi,\,\,\Psi
    \in\mathcal{F}^\ast_ {(\theta_1\theta_2)}
    (N_1'\times N_2').
\end{equation}\end{defn}
\begin{prop}\label{prop3}
The
 Gross Laplacian $\Delta_G$ on $\mathcal{F}^\ast_{(\theta_1,\theta_2)}(N_1'\times N_2')
 $ coincide with the adjoint of the
Gross Laplacian $\Delta_G^\ast$ on the test functions space
$\mathcal{F}_{(\theta_1,\theta_2)}(N_1'\times N_2')$.
\end{prop}
\textbf{Proof.}
For any $\Phi
\in\mathcal{F}^\ast_{(\theta_1,\theta_2)}(N_1'\times N_2')$ and
$\varphi\in\mathcal{F}_{(\theta_1,\theta_2)}(N_1'\times N_2')$, we
have
\begin{eqnarray*}
\langle\!\langle \Delta_G\Phi,\varphi\rangle\!\rangle
&=&\langle\!
\langle\mathbf{\mathcal{T}}\ast\Phi,\varphi\rangle\!\rangle  \\
   &=& [(\mathbf{\mathcal{T}}\ast\Phi)\ast\varphi](0)= [\Phi\ast(\mathbf{\mathcal{T}}\ast\varphi)](0)\\
   &=& \langle\!
\langle\Phi,\mathbf{\mathcal{T}}\ast\varphi\rangle\!\rangle=\langle\!
\langle\Phi,\Delta_G(\varphi)\rangle\!\rangle \\
&=&\langle\! \langle\Delta_G^\ast(\Phi),\varphi\rangle\!\rangle.
\end{eqnarray*}
\section{Quantum Gross Laplacian}\label{IV}
From the nuclearity of the space $\mathcal{F}_ {\theta_1}(N_1')$, we
have by Schwartz-Grothendieck Kernel Theorem the following
isomorphisms
\begin{equation}\label{isom}
    \mathcal{L}(\mathcal{F}_ {\theta_1}(N_1')
    ,\mathcal{F}^\ast_ {
\theta_2}(N_2')) \simeq\mathcal{F}^\ast_{\theta_1
}(N_1')\hat{\otimes}\mathcal{F}^\ast_{\theta_2 }(N_2')
 \simeq\mathcal{F}^\ast_{(\theta_1,
\theta_2)}(N_1'\times N_2') .
\end{equation}
So, for every $\Xi\in\mathcal{L}(\mathcal{F} _ {\theta_1}(N_1')
    ,\mathcal{F}^\ast_ {
\theta_2}(N_2'))$, the associated kernel $\Xi^K\in
\mathcal{F}^\ast_{(\theta_1, \theta_2)}(N_1'\times N_2')$ is defined
by
\begin{equation}\label{kernelprop}
    \langle\!\langle\Xi\varphi,\psi\rangle\!\rangle=
    \langle\!\langle\Xi^K,\varphi\otimes\psi\rangle\!\rangle
,\,\forall\varphi,\psi\in\mathcal{F}_{(\theta_1,
\theta_2)}(N_1'\times N_2').
\end{equation}
In the sequel, we will identify every operator $\Xi$ with its
kernel $\Xi^K.$

For every $\Xi\in\mathcal{L}(\mathcal{F} _ {\theta_1}(N_1')
    ,\mathcal{F}^\ast_ {
\theta_2}(N_2'))$, the symbol is defined in the usual manner (see
\cite{ob,o}\,) by
\begin{equation}\label{symboldef}
\sigma(\Xi)(\xi_1,\xi_2)=\langle\!\langle\Xi e_{ \xi_1}
,e_{\xi_2}\rangle\!\rangle=\langle\!\langle\Xi^K,e_
{(\xi_1,\xi_2)}\rangle\!\rangle=L(\Xi^K)(\xi_1,\xi_2)
\end{equation}
where $\xi_1\in N_1,\,\xi_2\in N_2$. Then, every operator in
$\mathcal{L}(\mathcal{F} _ {\theta_1}(N_1')
    ,\mathcal{F}^\ast_ {
\theta_2}(N_2'))$ is uniquely determined by its symbol since the
exponential vectors $e_{\xi_1}\otimes e_{\xi_2}$ span a dense
subspace of $\mathcal{F}_{(\theta_1, \theta_2)}(N_1'\times N_2')$.
As direct application of Theorem \ref{t3} and the identification
(\ref{kernelprop}) by the kernel theorem, we obtain the following
characterization of operators:
\begin{thm}\label{thm1}\cite{co}
The symbol map $\Xi\mapsto\sigma(\Xi)$ is a topological isomorphism:
\begin{eqnarray*}
\mathcal{L}(\mathcal{F}_ {\theta_1}(N_1')
    ,\mathcal{F}^\ast_ {
\theta_2}(N_2'))
&\rightarrow&\mathcal{G}_{(\theta^\ast_1,\theta^\ast_2)}(N_1\times
N_2) .
\end{eqnarray*}
\end{thm}

Let $\Xi_1$, $\Xi_2\in \mathcal{L}(\mathcal{F}_ {\theta} (N')
,\mathcal{F}_ {\theta} ( N'))$, we define the convolution product of
$\Xi_1$ and $\Xi_2$ denoted by $\Xi_1\ast \Xi_2$ by
\begin{equation}\label{convolutionoperator}
\sigma(\Xi_1\ast\Xi_2)=\sigma(\Xi_1)\sigma(\Xi_2).
\end{equation}
Using the topological isomorphisms:
\begin{equation}\label{isom1}
\mathcal{L}(\mathcal{F}_
{\theta_1} (N_1') ,\mathcal{F}^\ast_ {\theta_2} (
N_2'))\ni\Xi\mapsto\Xi^K\in\mathcal{F}^\ast_ {(\theta_1,\theta_2)}
(N_1'\times N_2')
\end{equation}
 defined in (\ref{isom}) and the extended Gross
Lapalcian given in (\ref{23}), we can define the quantum Gross
Laplacian as follows.

\begin{defn}
The quantum Gross
Laplacian $\Delta^Q_G$ is defined by:
\begin{equation}\label{40}
\Delta^Q_G(\Xi)=\mathbf{\mathcal{T}}\ast\Xi^K,\,\Xi\in\mathcal{L}(\mathcal{F}_ {\theta_1}(N_1')
    ,\mathcal{F}^\ast_ {
\theta_2}(N_2')),
\end{equation}where $\mathbf{\mathcal{T}}$ is defined by (\ref{24}).
\end{defn}
\begin{prop}\label{prop1}
For all $\Xi\in\mathcal{L}(\mathcal{F}_ {\theta_1}(N_1')
    ,\mathcal{F}^\ast_ {
\theta_2}(N_2'))$, the quantum Gross
Laplacian $\Delta^Q_G$ has the following analytic characterization:
\begin{equation}\label{symbolgross}
    \sigma(\Delta^Q_G(\Xi))(\xi_1,\xi_2)
    =(\langle\xi_1,\xi_1\rangle_1+\langle\xi_2
    ,\xi_2\rangle_2)\sigma(\Xi)(\xi_1,\xi_2),\,
    (\xi_1,\xi_2)\in N_1\times N_2.
\end{equation}
\end{prop}
\textbf{Proof.}
Using the property (\ref{convolutionoperator}) of the symbol of the convolution product
 of two operators and the definition (\ref{symboldef}), we have
\[\sigma(\Delta_G^Q(\Xi))=L(\mathcal{L}\ast\Xi^K)
=L(\mathbf{\mathcal{T}})L(\Xi^K)
=L(\mathbf{\mathcal{T}})\sigma(\Xi),\,\Xi\in\mathcal{L}(\mathcal{F}_ {\theta_1}(N_1')
    ,\mathcal{F}^\ast_ {
\theta_2}(N_2')).\]
 In the other hand, using the Taylor expansion of $\mathbf{\mathcal{T}}$, we obtain for all $(\xi_1
    ,\xi_2)\in N'_1\times N'_2$
    \begin{eqnarray*}
     L(\mathbf{\mathcal{T}}) (\xi_1
    ,\xi_2) &=&\langle\mathbf{\mathcal{T}}_{2,0}, \xi_1^{\otimes2}\rangle _1+ \langle\mathbf{\mathcal{T}}_{0,2}, \xi_2^{\otimes2}\rangle_2\nonumber\\
       &=&\langle\tau_{1}, \xi_1^{\otimes2}\rangle _1 +\langle\tau_{2}, \xi_2^{\otimes2}\rangle_2\nonumber\\&=&
       \langle\xi_{1}, \xi_1\rangle_1 +\langle\xi_{2}, \xi_2\rangle_2\nonumber.
    \end{eqnarray*}
    This proves the proposition.

\section{Relation between classical and quantum Gross Laplacian}\label{V}
In this section, we consider the particular case where
$\theta_1=\theta_2=\theta$ and $N=N_1=N_2$. For simplicity, we use
the same notation for the Gross Laplacian acting on
$\mathcal{F}_{\theta}(N')$ and on
$\mathcal{F}_{(\theta_1,\theta_2)}(N_1'\times N_2')$. It is well
known (see \cite{ghor}), that the pointwise multiplication yields a
continuous bilinear map from
$\mathcal{F}_{\theta}(N')\times\mathcal{F}_{\theta}(N')$ into
$\mathcal{F}_{\theta}(N')$, i. e., for all
$f,g\in\mathcal{F}_{\theta}(N')$, we have
$fg\in\mathcal{F}_{\theta}(N')$. Let
$\Phi\in\mathcal{F}_{\theta}^\ast(N')$ fixed and
$\mathcal{M}_\Phi\in\mathcal{L}(\mathcal{F}_{\theta}
(N'),\mathcal{F}_{\theta}^\ast(N'))$ be the multiplication operator
by $\Phi$ defined by
\begin{equation}\label{defmulti}
\langle\!\langle\mathcal{M}_\Phi f,g\rangle\!\rangle=
\langle\!\langle \Phi,fg\rangle\!\rangle,\,\,f,g\in
\mathcal{F}_{\theta}(N').
\end{equation}
Let
$e_0=(1,0,\ldots)\in\mathcal{F}_{\theta}(N')$ be the vacuum vector. Using the definition (\ref{defmulti}), we have the following connection between the multiplication operator and a distribution given by
$$\mathcal{M}_\Phi e_0=\Phi.$$Therefore, we obtain the following relation between classical and quantum Gross Laplacian.
\begin{prop}
For any $\Phi\in\mathcal{F}_{\theta}^\ast(N')$, we have:
\begin{equation}\label{47}
    \left(\Delta^Q_G\mathcal{M}_{\Phi}\right)e_0=
    \Delta_G\Phi.
\end{equation}
\end{prop}
\textbf{Proof.} Let $\Phi\in\mathcal{F}_{\theta}^\ast(N')$. Using the definitions of the Laplace transform (\ref{2'}) and the multiplication operator (\ref{defmulti}), we have
\begin{eqnarray}
   \sigma(\mathcal{M}_\Phi)(\xi,\eta)&=&
   \langle\!\langle
\mathcal{M}_\Phi e_\xi, e_\eta\rangle\!\rangle  \nonumber\\
   &=&\langle\!\langle
\Phi,e_{\xi+\eta}\rangle\!\rangle  \nonumber\\
   &=& \mathrm{L} (\Phi
)(\xi+\eta),\,\xi,\eta\in N\label{laplacemultiplication}.
\end{eqnarray}In the other hand, the Laplace transform of the classical Gross Lapalcian is given by
$$L(\Delta_G(\Phi))(\eta)=L(\mathcal{T}\ast\Phi)(\eta)=
\langle\eta,\eta\rangle\widehat{\Phi}(\eta),\,\eta\in N.$$
By Proposition \ref{prop1} and the equality (\ref{laplacemultiplication}), we obtain
\begin{eqnarray*}
L(\Delta_G^Q(\mathcal{M}_\Phi)(e_0))(\eta)&=&\langle\langle\Delta_G^Q(\mathcal{M}_\Phi)(e_0),e_{\eta}
\rangle \rangle  \nonumber\\ &=&\langle\langle\Delta_G^Q\mathcal{M}_\Phi,e_0\otimes
e_{\eta}
\rangle\rangle\nonumber  \\
   &=& \sigma(\Delta_G^Q(\mathcal{M}_\Phi))(0,\eta)\nonumber \\
   &=&\langle\eta,\eta\rangle\sigma (\mathcal{M}_\Phi))(0,\eta)\nonumber \\
   &=&\langle\eta,\eta\rangle\widehat{\Phi}(\eta) \nonumber\\&=&L(\Delta_G(\Phi))(\eta),\,\eta\in N\nonumber.
\end{eqnarray*}
Using the fact that the Laplace transform is an isomorphism, we have
\[\Delta_G^Q\mathcal{M}_\Phi e_0=\Delta
_G\Phi,\,\Phi\in \mathcal{F}_{\theta}^\ast(N').\]

\section{Solution of linear quantum
stochastic differential equations}\label{VI} Let
$I\subset\mathbb{R}$ be an interval containing the origin. Consider
a family $\{\Phi_t;t\in I\}$ of distributions in
$\mathcal{F}^\ast_{(\theta_1,\theta_2)}(N_1'\times N_2')$ i. e., the
function $t\mapsto\Phi_t$ is continuous from $I$ into
$\mathcal{F}^\ast_{(\theta_1,\theta_2)}(N_1'\times N_2')$. So, by
using the isomorphism between
$\mathcal{F}^\ast_{(\theta_1,\theta_2)}(N_1'\times N_2')$ and
$\mathcal{G}_{(\theta^\ast_1,\theta^\ast_2)}(N_1\times N_2)$ via the
Laplace transform, the function $t\mapsto \widehat{\Phi}_t$ is
continuous from $I$ into
$\mathcal{G}_{(\theta^\ast_1,\theta^\ast_2)}(N_1\times N_2)$. Then,
for each $t\in I$, the set $\{\widehat{\Phi}_s;s\in [0,t]\}$ is a
compact subset of the generalized space
$\mathcal{G}_{(\theta^\ast_1,\theta^\ast_2)}(N_1\times N_2)$. This
implies that it is bounded in
$\mathcal{G}_{(\theta^\ast_1,\theta^\ast_2)}(N_1\times N_2)$. Hence,
there exist constants $p,q\in\mathbb{N}$, $a_1,a_2>0$ and $C_t>0$
such that
\[|\widehat{\Phi}_t(\xi,\eta)|\leq C_te^{\theta_1^\ast
(a_1 |\xi|_{p})+\theta_2^\ast(a_2 |\eta|_{q})}|,\,(\xi,\eta)\in
N_1\times N_2.\]By consequence, the function
$(\xi,\eta)\mapsto\int_0^t\Phi_s(\xi,\eta)ds$ and belongs to the space
$\mathcal{G}_{(\theta^\ast_1,\theta^\ast_2)}(N_1\times N_2)$. This
way we can define the integral $E_t=\int_0^t\Phi_s(\xi,\eta)ds$ as
the unique element of $\mathcal{F}^
\ast_{(\theta_1,\theta_2)}(N_1'\times N_2')$ satisfying
\[\left(\int_0^t\Phi_s(\xi,\eta)ds\right)^{\widehat{}}
=\int_0^t\widehat{\Phi}_s(\xi,\eta)ds,\,(\xi,\eta)\in N_1\times
N_2.\] Moreover, for $t\in I$, the process $E_t$ is differentiable
in $\mathcal{F}_ {\theta}(N'\times N')^ \ast$ and satisfies the
equation
\[\frac{\partial E_t}{\partial t}=\Phi_t.\]

\begin{prop}\label{l1}
For every distribution
$\Phi\in\mathcal{F}^ \ast_ {(\theta_1,\theta_2)}(N_1'\times N_2')$,
the functional $e^{\ast\Phi}$ is defined by
\begin{equation}\label{25}
\widehat{e^{\ast\Phi}}=e^{\widehat{\Phi}}
\end{equation}belongs to $\mathcal{F}^ \ast_ {\left((e^{\theta_1^\ast})^\ast,(e^{\theta_2^\ast})^\ast
\right)}(N_1'\times N_2').$\end{prop} \begin{proof}The proof is
similar to the one of Theorem 1 established in the paper \cite{coo}.
\end{proof} Consider the following initial value problem:
\begin{equation}\label{43}
\left\{
  \begin{array}{ll}
 \frac{d\Xi(t)}{dt}=Z(t)\ast\Xi(t)+\Theta(t)   \\
  \Xi(0)=\Xi_0,
  \end{array}
\right.
\end{equation}
where $t\mapsto Z(t)$ and $t\mapsto\Theta(t)$ are continuous map defined on an interval
$I$ into $\mathcal{L}(\mathcal{F}_ {\theta_1}(N_1')
    ,\mathcal{F}^\ast_ {
\theta_2}(N_2'))$ and $\Xi_0\in\mathcal{L}(\mathcal{F}_ {\theta_1}(N_1')
    ,\mathcal{F}^\ast_ {
\theta_2}(N_2')).$

\begin{thm}\label{t5}The stochastic quantum differential equation
(\ref{43}) has a unique solution in $\mathcal{L}(\mathcal{F} _
{\left(e^{\theta_1^\ast}\right)^\ast}(N_1')
    ,\mathcal{F}^\ast_ {
\left(e^{\theta_2^\ast}\right)^\ast}(N_2'))$ given
by\begin{equation}\label{44} \Xi(t)=\Xi_0\ast e^{\ast\int_0^tZ(s)ds}
    +\int_0^te^{\ast( \int_s^tZ(u)du)}
    \ast\Theta^K_sds.
\end{equation}
\end{thm}
\textbf{Proof.} Applying the symbol map to the differential equation (\ref{43}), we obtain an ordinary differential equation given by
\begin{equation}\label{44}
\left\{
  \begin{array}{ll}
 \frac{d\sigma(\Xi)(t)}{dt}=\sigma(Z)(t)\sigma(\Xi)(t)+\sigma(
 \Theta)(t)   \\
  \sigma(\Xi)(0)=\sigma(\Xi_0)\in\mathcal{F}^\ast_{(\theta_1,\theta_2)}(N'_1\times N'_2),
  \end{array}
\right.
\end{equation}
whose solution is given by
$$\sigma(\Xi)(t)=\sigma(\Xi_0) e^{\int_0^{t}\sigma(Z)(s)ds}
    +\int_0^{t}e^{\int_s^t\sigma(Z)(u)du}
    \sigma(\Theta)(s)ds.$$Moreover, by Proposition \ref{l1}, $\sigma(\Xi)(t)$ is an element of $\mathcal{F}^
\ast_ {\left(\left(e^{\theta_1^\ast}\right)^\ast,\left(e^{\theta_2
^\ast}\right) ^\ast\right)} (N_1'\times N_2')$, for all $t\in I$.
Using the property (\ref{convolutionoperator}), the solution of the linear quantum stochastic differential equation is given by
$$\Xi(t)=\Xi_0\ast e^{\ast\int_0^tZ(s)ds}
    +\int_0^te^{\ast( \int_s^tZ(u)du)}
    \ast\Theta_sds.$$
Now, by Theorem \ref{thm1}, we obtain that $$\Xi(t)\in\mathcal{L}(\mathcal{F}_ {\left(e^{\theta_1^\ast}\right)^\ast} (N_1'),\mathcal{F}^
\ast_ {\left(e^{\theta_2
^\ast}\right) ^\ast} ( N_2')),$$ for all $t\in I$.
\subsection{Solution of quantum Gross heat equation}
In this section, we consider two Young functions $\theta_1$ and $\theta_2$, satisfying
\[\displaystyle\lim_{x\rightarrow\infty}
\frac{\theta_i(x)}{x^2}<\infty,\,i=1,2.\]

\begin{thm}\label{t6}The quantum Gross heat equation
\begin{equation}\label{45}
\left\{
  \begin{array}{ll}
 \frac{d\Xi(t)}{dt}=\frac{1}{2}\Delta_G^Q\Xi(t)+\Theta(t)   \\
  \Xi(0)=\Xi_0,
  \end{array}
\right.
\end{equation} has a unique solution in $\mathcal{L}(\mathcal{F} _
{\left(e^{\theta_1^\ast}\right)^\ast}(N_1')
    ,\mathcal{F}^\ast_ {
\left(e^{\theta_2^\ast}\right)^\ast}(N_2'))$ given
by\begin{equation}\label{44} \Xi(t)=\Xi_0\ast e^{\ast\frac{t}{2}\mathbf{\mathcal{T}}}
    +\int_0^te^{\ast\frac{t-s}{2}\mathbf{\mathcal{T}}}
    \ast\Theta(s)ds.
\end{equation}
\end{thm}
\begin{proof}
For all $t\in I,$ let $Z(t)=\mathbf{\mathcal{T}}$, where $\mathbf{\mathcal{T}}$ is the distribution given by (\ref{24}). Then, the equation (\ref{44}) becomes the quantum Gross heat equation given by
\begin{equation}\label{45'}
\left\{
  \begin{array}{ll}
 \frac{d\Xi(t)}{dt}=
 \frac{1}{2}\mathbf{\mathcal{T}}\ast\Xi(t)
+\Theta(t) =\frac{1}{2}\Delta_G^Q\Xi(t)+\Theta(t)   \\
  \Xi(0)=\Xi_0.
  \end{array}
\right.
\end{equation} Therefore, we apply Theorem \ref{t5} to get the unique solution in equation (\ref{44}). We can further rewrite this solution in another way. For $t>0$, we define $\gamma_t=\gamma_1^t\otimes\gamma_2^t(.)=\gamma_1\otimes
\gamma_2(./\sqrt{t})$, where $\gamma=\gamma_1\otimes\gamma_2$ is the standard Gaussian measure on the space $X'$ defined in (\ref{gaussian measure}). It follows that the solution (\ref{44}) can be expressed as \[\Xi(t)=\Xi_0\ast \widetilde{\gamma}_t
    +\int_0^t\widetilde{\gamma}_{t-s}
    \ast\Theta(s)ds,\]where $\widetilde{\gamma}$ is a positive distribution in $\mathcal{F}^\ast_ {(\theta_1,\theta_2)} (N_1'\times N'_2)$ given by
    \[\langle\!\langle\widetilde{\gamma}_{t},
    \varphi\rangle\!\rangle=\int_{X'_1\times X'_2}\varphi(x_1,x_2)d\gamma_t(x_1,x_2)=\int_{X'_1\times X'_2}\varphi(\sqrt{t}x_1,\sqrt{t}x_2)d\gamma_1(x_1)d\gamma_2(x_2),\]
where $\varphi\in\mathcal{F}_ {(\theta_1,\theta_2)} (N_1'\times N'_2).$
\end{proof}It is easy to see from the definitions of the kernel operator (\ref{kernelprop}) and the convolution product of the operators (\ref{convolutionoperator}), that for all $\Xi_1,\Xi_2\in\mathcal{L}(\mathcal{F} _
{\theta_1}(N_1')
    ,\mathcal{F}^\ast_ {
\theta_2}(N_2'))$, we have
\[\Xi_1\ast\Xi_2=\Xi_1^K\ast\Xi_2^K.\]
Therefore combining Theorem \ref{t6} for the particular case where $\Theta=0$ and the Gel'fand triplet (\ref{triplet2}), we obtain the following result:
\begin{cor}Let $\Xi_0\in\mathcal{F} _
{(\theta_1,\theta_2)}(N_1'\times N_2')$. Then,
the quantum Gross heat equation
\begin{equation}\label{45}
\left\{
  \begin{array}{ll}
 \frac{d\Xi_t}{dt}=\frac{1}{2}\Delta_G^Q\Xi(t) \\
  \Xi(0)=\Xi_0,
  \end{array}
\right.
\end{equation} has a unique solution whose kernel is given
by\begin{equation}\label{44} \Xi^K_t(y_1,y_2)=\int_{X_1\times X_2}
\Xi_0(y_1+\sqrt{t}x_1,y_2+\sqrt{t}x_2)
d\gamma_1(x_1)d\gamma_2(x_2)
.
\end{equation}
\end{cor}


\begin{thebibliography}{99}
\addcontentsline{toc}{chapter}{\quad\ Bibliographie}
\bibitem{bko}\textbf{A. Barhoumi, H. H. Kuo and H. Ouerdiane},
\emph{Generalized Gross heat equation with noises}, Soochow J.
Math., V. 32, No. 1, (2006), 113-125.
\bibitem{CJ99}\textbf{D. M. Chung and U. C. Ji}, \emph{Some Cauchy problems in white noise
analysis and associated semigroups of operators}, Stochastic Anal.
Appl., 17(1), (1999), 1-22.
\bibitem{co}\textbf{M. Ben Chrouda and H. Ouerdiane}, \emph{Algebra of operators on holomorphic
functions and Applications}, Journal of Mathematical Physics,
Analysis and Geometry, Vol. 5 (2002), 65-76.
\bibitem{coo}\textbf{M. Ben Chrouda, M. Eloued and H. Ouerdiane},
\emph{Convolution calculus and applications to stochastic
differential equations}, Soochow J. Math., 28(4), (2002), 375-388.
\bibitem{ghor}\textbf{R. Gannoun, R. Hachaichi, H. Ouerdiane and A. Rezgui
}, \emph{Un th\'{e}or\`{e}me de dualit\'{e} entre espaces de
fonctions holomorphes \`{a} croissance exponentielle }, J. Func.
Anal. vol. 171, no. 1, (2000), 1-14.
\bibitem{g}\textbf{L. Gross}, \emph{Potential theory on Hilbert space}, J. Funct. Anal. , 1, (1967), 123-
181.

\bibitem{HKPS93}\textbf{T. Hida, H. H. Kuo, J. Potthoff, and L. Streit}, \emph{White Noise. An
Infinite Dimensional Calculus}, Kluwer Academic Publishers,
Dordrecht, (1993).
\bibitem{HOS92}\textbf{T. Hida, N. Obata, and K. Saitô}, \emph{Infinite-dimensional rotations
and Laplacians in terms of white noise calculus}, Nagoya Math. J. ,
128, (1992), 65-93.
\bibitem{joo}\textbf{U. C. Ji, N. Obata and H. Ouerdiane}, \emph{Analytic
Characterization of Generalized Fock Space Operators as Two-variable
Entire Functions with Growth Condition}, World Scientific, vol. 5,
No. 3, (2002), 395-407.
\bibitem{k'}\textbf{ H. H. Kuo},
\emph{On Laplacian operators of generalized Brownian functionals},
Lect. Notes in Math. 1203 (1986), 119-128.
\bibitem{k}\textbf{ H. H. Kuo},
\emph{White noise distribution theory}, CRC Press, (1996).
\bibitem{ob}\textbf{N. Obata}, \emph{White Noise Calculus and Fock
Spaces}, Lecture Notes in Math. 1577 (1994).
\bibitem{o}\textbf{H. Ouerdiane}, \emph{Noyaux et symboles d'op\'{e}rateurs sur des fonctionnelles analytiques
gaussiennes}, Japan. J. Math. 21 (1995), 223-234.
\bibitem{ou}\textbf{H. Ouerdiane}, \emph{Infinite dimensional
entire functions and applications to stochastic differential
equations}, Notices of th South African Math. Society, 35, No. 1
(2004), 23-45.
\bibitem{pie70}\textbf{M. A. Piech}, \emph{A fundamental solution of the parabolic equation on
Hilbert space II. The semigroup property}, Trans. Amer. Math. Soc.
150, (1970), 257-286.
\end{thebibliography}
\end{document}